\documentclass[leqno,11pt]{article}
\usepackage{amssymb,latexsym}
\usepackage[dvips]{graphicx} 
\usepackage{url}

 \catcode`@=11\def\@oddhead{\hbox{}\hfil\rm\thepage}\def\@oddfoot{}
 \def\@evenhead{\hbox{}\hfil\rm\thepage}\def\@evenfoot{}
 \@twosidetrue \catcode`@=12
\newtheorem{prp}{Proposition}
\newtheorem{lem}[prp]{Lemma}\newtheorem{thm}[prp]{Theorem}
\newtheorem{cor}[prp]{Corollary}
\newenvironment{proof}{\begin{trivlist}\item[\textit{Proof.}]}{\end{trivlist}
  \medskip\par}
\newenvironment{prfof}[1]{\begin{trivlist}\item[\textit{Proof of #1.}]}{
  \end{trivlist} \medskip \par}
\newenvironment{rem}{\begin{trivlist}\item[\textit{Remark.}]}{\end{trivlist}
  \medskip\par}
\newenvironment{xmp}{\begin{trivlist}\item[\emph{Example. }]}{\end{trivlist}
  \bigskip\par}
\def\prpb{\begin{prp}}\def\prpe{\end{prp}}
\def\lemb{\begin{lem}}\def\leme{\end{lem}}
\def\thmb{\begin{thm}}\def\thme{\end{thm}}
\def\corb{\begin{cor}}\def\core{\end{cor}}
\def\xmpb{\begin{xmp}}\def\xmpe{\end{xmp}}
\def\prfb{\begin{proof}}\def\prfe{\end{proof}}
\def\prfofb#1{\begin{prfof}{#1}}\def\prfofe{\end{prfof}}
\def\remb{\begin{rem}}\def\reme{\end{rem}}
\def\prpa#1{\label{p:#1}}\def\prpu#1{Proposition~\ref{p:#1}}
\def\lema#1{\label{l:#1}}\def\lemu#1{Lemma~\ref{l:#1}}
\def\thma#1{\label{t:#1}}\def\thmu#1{Theorem~\ref{t:#1}}
\def\cora#1{\label{c:#1}}\def\coru#1{Corollary~\ref{c:#1}}
\def\seca#1{\label{s:#1}}\def\secu#1{\S~\ref{s:#1}}

\def\itmb{\begin{enumerate}}\def\itme{\end{enumerate}}
\def\itdb{\begin{itemize}}\def\itde{\end{itemize}}
\def\ittb{\begin{description}}\def\itte{\end{description}}
\def\arrb#1{\begin{array}{#1}}\def\arre{\end{array}}
\def\tabb#1{\par\noindent\begin{tabular}{#1}}
\def\tabe{\end{tabular}\par\noindent}
\def\eqna#1{\label{e:#1}}\def\eqnu#1{(\ref{e:#1})}

\def\QED{\relax\ifmmode\let\@tempa\relax\ifcase\@eqcnt\def\@tempa{& & &}\or
  \def\@tempa{& &}\else\def\@tempa{&}\fi\@tempa $\Box$ \else\hfill $\Box$ \fi}
\def\DDD{\relax\ifmmode\let\@tempa\relax\ifcase\@eqcnt\def\@tempa{& & &}\or
 \def\@tempa{& &}\else\def\@tempa{&}\fi\@tempa $\Diamond$
 \else\hfill $\Diamond$ \fi}

\def\Rom#1{\uppercase\expandafter{\romannumeral#1}}
\def\dsp{\displaystyle}
\def\eps{\epsilon}

\def\limf#1{\displaystyle \lim_{#1\to\infty}}

\def\Ccomb#1#2{\setbox0=\hbox{$\displaystyle\mathrm{C}$}\setbox1=\hbox{%
$\scriptstyle #1$}\kern \wd1{\mathrm{C}}_{\kern -1.05\wd0\kern -0.99\wd1{#1}
 \kern 1.15\wd0{#2}}}

\def\clvec#1#2#3{\def\clvecone{#3}\left(\arrb{c} \dsp #1\\ \dsp #2
 \ifx\clvecone\empty\else\\ \dsp #3\fi\arre\right)}

\def\le{\leqq} \def\ge{\geqq} 

\def\reals{{\mathbb R}}

\def\nintegers{{\mathbb N}}

\def\prb#1{\def\prbone{#1}
  \ifx\prbone\empty{\mathrm{P}}\else{\mathrm{P[\;}}#1{\mathrm{\;]}}\fi}
\def\prbseq#1#2{\def\prbseqone{#2}
  \ifx\prbseqone\empty{\mathrm{P}}_{#1}\ignorespaces
  \else{\mathrm{P}}_{#1}{\mathrm{[\;}}#2{\mathrm{\;]}}\fi}

\def\EEseq#1#2{\def\EEseqone{#2}
  \ifx\EEseqone\empty{\mathrm{E}}_{#1}\else
 {\mathrm{E}}_{#1}{\dsp\mathrm{[\;}}#2{\mathrm{\;]}}\fi}

\def\VVseq#1#2{\def\VVseqone{#2}
  \ifx\VVseqone\empty{\matrm{V}}_{#1}\else
 {\mathrm{V}}_{#1}{\dsp\mathrm{[\;}}#2{\mathrm{\;]}}\fi}

\def\chrfcn#1{\mathop{\mathbf{1}}\nolimits_{#1}}

\def\subsetne{\mbox{\raisebox{0.15em}{\setbox0=\hbox{$\subset$}\copy0\kern-0.5\wd0}\raisebox{-0.60em}{\setbox0=\hbox{$\ne$}\kern-0.5\wd0\copy0\kern-0.5\wd0}\setbox0=\hbox{$\subset$}\kern+0.5\wd0}}

\title{An elementary proof of representation of submodular function as an supremum of measures on $\sigma$-algebra with totally ordered generating class}
\author{
Tetsuya Hattori
\footnote{ 
This work was supported by JSPS KAKENHI Grant Number 22K03358.
}
\small Keio University,
\\ \small URL: \url{https://tetshattori.web.fc2.com/research.htm}
\\ \small email: \url{hattori@econ.keio.ac.jp}
}
\date{2023/10/21}

\begin{document}

\maketitle

\begin{abstract}

We give an alternative proof of a fact that a 
finite continuous non-decreasing submodular set function 
on a measurable space
can be expressed as a supremum of measures dominated by the function,
if there exists a class of sets which is totally ordered with respect to
inclusion and generates the sigma-algebra of the space.
The proof is elementary in the sense that 
the measure attaining the supremum in the claim is 
constructed by a standard extension theorem of measures.
As a consequence, a uniquness of the supremum attaining measure also follows.
A Polish space is an examples of the measurable space
which has a class of totally ordered sets that 
generates the Borel sigma-algebra.

\end{abstract}


2020 MSC Primary 60A10; Secondary 60Axx

\section{Introduction}

Let $(\Omega,{\cal F})$ be a measurable space,
namely, a $\sigma$-algebra ${\cal F}$ is a class of
subsets of $\Omega$ and is closed under 
complements and countable unions.
For a measurable set $A\in{\cal F}$ denote by ${\cal F}|_{A}$,
the class of measurable sets restricted to $A$, and
denote the set of finite measures on the measurable space $(A,{\cal F}|_A)$
by ${\cal M}(A)$.

For a set function $v:\ {\cal F}\to\reals$ and
a measurable set $A\in{\cal F}$,
let ${\cal C}_{-,v}(A)$ be a class of measures defined by
\begin{equation}
\eqna{submodularlocalcoredef}
{\cal C}_{-,v}(A)=\{\mu\in{\cal M}(A)
\mid \mu(A)=v(A),
\ \mu(B)\le v(B),\ B\in {\cal F}|_{A}\}.
\end{equation}
If
\begin{equation}
\eqna{submodulartotalcorecoherentriskdef}
v(B)=\sup_{\mu\in{\cal C}_{-,v}(A)} \mu(B),
\end{equation}
holds for all $A,B\in{\cal F}$ satisfying $B\subset A$,
then it is easy to see that
\begin{equation}
\eqna{submodular}
v(A)+v(B)\ge v(A\cup B)+v(A\cap B), \ \ A,B\in {\cal F},
\end{equation}
holds (see \prpu{totalcoherentriskismodular} in \secu{defandresult}).
A set function which satisfies \eqnu{submodular} is called
a submodular function.

The converse that a submodular function satisfies
\eqnu{submodulartotalcorecoherentriskdef} is also known to hold
under mild and natural assumptions.
In fact, a proof in \cite{Denneberg} proves existence of a measure $\mu$
which satisfies $v(B)=\mu(B)$, so that the supremum in
\eqnu{submodulartotalcorecoherentriskdef} is attained.

The formula \eqnu{submodulartotalcorecoherentriskdef}
has significance in the related fields of study, such as 
coherent risk measures in mathematical finance and
cores of convex games in cooperative game theory.
In view of a wide interest in this formula,
it may be worthwhile to find an alternative elementary proof.

In contrast to a proof in \cite{Denneberg} which 
seeks for wide applicability even beyond measurable spaces,
we keep ourselves as close as possible to measures,
except for \eqnu{submodular} which characterizes the submodular property.
See the definitions in \secu{defandresult} for details.
Our proof in \secu{proof} is elementary in the sense that we
prove the existence of $\mu$ satisfying $v(B)=\mu(B)$ by
the extension theorem of measures on a finite algebra to 
the $\sigma$-algebra generated by the finite algebra,
in contrast to the proof in \cite{Denneberg} which
uses Hahn--Banach's Theorem. As a consequence, uniqueness of $\mu$ 
do not follow in general in the latter proof,
while we have certain uniqueness result for $\mu$
(see \thmu{totorderimpliesgenDelbaenKusuoka}
in \secu{defandresult}).
In this uniqueness result, it is essential that
there exists of a class 
${\cal I}$ which is totally ordered with respect to inclusion
and generates the $\sigma$-algebra ${\cal F}$.
Examples of measurable spaces $\dsp(\Omega,{\cal F})$ which has
such a class ${\cal I}$ are given in \secu{salggeneratingchain}.
Polish spaces are in the examples, hence our main result holds for spaces 
which are extensively used in the theory of stochastic processes.

\section{Definition and main result}
\seca{defandresult}

Given a set function $v:\ {\cal F}\to\reals$ 
we define to say that $v$ is non-decreasing, if $v$ satisfies
\begin{equation}
\eqna{nondecreasing}
v(A)\le v(B),\ A\subset B,\ A,B\in{\cal F},
\end{equation}
holds, and that $v$ is continuous, if
\begin{equation}
\eqna{upperconti}
\limf{n} v(A_n)=v(\bigcup_{n\in\nintegers} A_n),
\ \ A_1\subset A_2\subset\cdots,\ A_n\in{\cal F},\ n=1,2,3,\ldots,
\end{equation}
and 
\begin{equation}
\eqna{lowerconti}
\limf{n} v(A_n)=v(\bigcap_{n\in\nintegers} A_n)
\ \ A_1\supset A_2\supset\cdots,\ A_n\in{\cal F},\ n=1,2,3,\ldots,
\end{equation}
hold.

In measure theory,
a finite measure is defined to be a non-negative real valued
$\sigma$-additive set function defined on a $\sigma$-algebra ${\cal F}$.
Equivalently, we can say that
a real valued set function $\mu:\ {\cal F}\to\reals$
is a finite measure, 
if $\mu$ is non-decreasing, continuous, $\mu(\emptyset)=0$, and
\begin{equation}
\eqna{modular}
\mu(A)+\mu(B)= \mu(A\cup B)+\mu(A\cap B), \ \ A,B\in {\cal F}.
\end{equation}
Submodular and supermodular functions are defined by replacing 
the equality \eqnu{modular} with inequalities.
In this paper, 
we say that a set function $v:\ {\cal F}\to\reals$ is submodular,
if $v$ is non-decreasing, continuous, $v(\emptyset)=0$, and satisfies
\eqnu{submodular},
and $v:\ {\cal F}\to\reals$ is supermodular,
if $v$ is non-decreasing, continuous, $v(\emptyset)=0$, and satisfies
\begin{equation}
\eqna{supermodular}
v(A)+v(B)\le v(A\cup B)+v(A\cap B), \ \ A,B\in {\cal F}.
\end{equation}

Note that, to keep the definitions close to that of measures,
we assume non-decreasing property \eqnu{nondecreasing} and
continuity \eqnu{upperconti} and \eqnu{lowerconti}
in the definitions of submodular and supermodular functions.
Note also that while $\mu(\emptyset)=0$ is essential for a measure $\mu$ to
have additivity property, submodular and supermodular functions lack
additivity property in general. Imposing $v(\emptyset)=0$ is thus
only for notational simplicity, and 
the formula in this paper can easily be generalized to the case 
$v(\emptyset)\ne 0$ by the replacements $v(A)\mapsto v(A)-v(\emptyset)$.

The results in this paper for supermodular functions are proved by
correspoinding results for submodular functions if we note the following
well-known relation.
\prpb
\prpa{convexconcaveduality}
For a set functoin $v:\ {\cal F}\to\reals$ on a measurable space 
$(\Omega,{\cal F})$ define a set function $\tilde{v}:\ {\cal F}\to\reals$ by
\begin{equation}
\eqna{convexconcavedualtr}
\tilde{v}(A)=v(\Omega)-v(A^c)+v(\emptyset),\ A\in{\cal F}.
\end{equation}
Then the following hold:
\itmb
\item
$\tilde{v}(\emptyset)=v(\emptyset)$ and $\tilde{v}(\Omega)=v(\Omega)$.
\item
If $v$ is non-decreasing, $\tilde{v}$ is also non-decreasing.
\item
If $v$ is continuous, $\tilde{v}$ is also continuous.
\item
If $v$ is submodular then $\tilde{v}$ is supermodular,
and conversly
if $v$ is supermodular then $\tilde{v}$ is submodular.
\DDD
\itme
\prpe
\prfb
The first three properties are obvious.
The last property is also easy, if we note that
if $v$ satisfies \eqnu{submodular} (or respectively, \eqnu{supermodular}),
the direction of the inequality for $\tilde{v}$ changes 
because of the minus sign in \eqnu{convexconcavedualtr}, together with
the de Morgan's laws and that ${\cal F}$ is closed under complement.
\QED
\prfe
In analogy to $\dsp{\cal C}_{-,v}(A)$ in \eqnu{submodularlocalcoredef},
let
\begin{equation}
\eqna{restriction}
{\cal F}|_A=\{B\in{\cal F}\mid B\subset A\}=\{B\cap A\mid B\in {\cal F} \}
\end{equation}
and let
${\cal C}_{+,v}(A)$ be a class of measures defined by
\begin{equation}
\eqna{localcoredef}
{\cal C}_{+,v}(A)=\{\mu\in{\cal M}(A)
\mid \mu(A)=v(A),
\ \mu(B)\ge v(B),\ B\in {\cal F}|_{A}\}.
\end{equation}
\prpb
\prpa{totalcoherentriskismodular}
Let $(\Omega,{\cal F})$ be a measureable space,
and $v:\ {\cal F}\to\reals$ a set function.
Then the following hold.
\itmb
\item
If \eqnu{submodulartotalcorecoherentriskdef}
holds for all $A,B\in{\cal F}$ satisfying $B\subset A$,
then $v$ is a submodular function.

\item
If
\begin{equation}
\eqna{totalcorecoherentriskdef}
v(B)=\inf_{\mu\in{\cal C}_{+,v}(A)} \mu(B)
\end{equation}
holds for all $A,B\in{\cal F}$ satisfying $B\subset A$,
then $v$ is a supermodular function.
\DDD
\itme
\prpe
%
%
\prfb
Assume \eqnu{submodulartotalcorecoherentriskdef}.
Since a measure $\mu$ is non-decreasing and continuous, 
$\dsp v(B)=\sup_{\mu\in{\cal C}_{-,v}(\Omega)} \mu(B)$ obtained
by $A=\Omega$ in \eqnu{submodulartotalcorecoherentriskdef}
implies that $v$ is also non-decreasing and continuous.

To prove \eqnu{submodular}, let $A,B\in {\cal F}$ and 
substitute $A$ in \eqnu{submodulartotalcorecoherentriskdef} by
$A\cup B$ to obtain
$\dsp v(B)=\sup_{\mu\in{\cal C}_{-,v}(A\cup B)} \mu(B)$
and $\dsp v(A)=\sup_{\mu\in{\cal C}_{-,v}(A\cup B)} \mu(A)$,
which further imply $v(B)\ge \mu(B)$ and $v(A)\ge \mu(A)$ for 
all $\mu\in {\cal C}_{-,v}(A\cup B)$.
Also \eqnu{submodularlocalcoredef} implies
$v(A\cup B)=\mu(A\cup B)$ for all $\mu\in {\cal C}_{-,v}(A\cup B)$.
Finally, for $\eps>0$,
$\dsp v(A\cap B)=\sup_{\mu\in{\cal C}_{-,v}(A\cup B)} \mu(A\cap B)$
implies that there exists $\mu\in{\cal C}_{-,v}(A\cup B)$ such that
$\dsp v(A\cap B)\le \mu(A\cap B)+\eps$.
These equality and inequalities imply, with \eqnu{modular},
\[
v(A\cup B)+v(A\cap B)-v(A)-v(B) \le \mu(A\cup B)+\mu(A\cap B)+\eps-\mu(A)-\mu(B)
=\eps.
\]
$\eps$ can be any positive constant, hence \eqnu{submodular} follows.

The second part of the Proposition follows from the first part with
\prpu{convexconcaveduality}.
\QED
\prfe

Concerning the converse of \prpu{totalcoherentriskismodular},
note first that the definition \eqnu{submodularlocalcoredef} of
${\cal C}_{-,v}(A)$ implies 
$\dsp v(B)\ge \sup_{\mu\in{\cal C}_{-,v}(A)} \mu(B)$
for any $A,B\in{\cal F}$ satisfying $B\subset A$.
Hence if there exists $\dsp\mu\in{\cal C}_{-,v}(A)$ such that 
$\mu(B)=v(B)$ then 
\eqnu{submodulartotalcorecoherentriskdef} holds for this pair $(A,B)$.
The main result of this paper is on the construction of such $\mu$.
To state the main theorem, we consider the following set of conditions on
a class of measurable sets ${\cal I}\subset{\cal F}$;
\begin{equation}
\eqna{sigmaalgebrawithtotallyorderedgeneratingclass}
\left\{\arrb{l}\dsp
\sigma[{\cal I}]={\cal F},
\mbox{ where $\sigma[{\cal I}]$ denotes the smallest $\sigma$-algebra
 containing ${\cal I}$,}
\\ \dsp
\emptyset\in{\cal I},\ \Omega\in{\cal I},\\ \dsp
I\mbox{ is totally ordered with respect to inclusion,}
\\ \dsp \mbox{i.e., for all $I_1,I_2\in{\cal I}$
either $I_1\subset I_2$ or $I_2\subset I_1$\,,}
\arre\right.
\end{equation}
%

For a class ${\cal I}\subset {\cal F}|_A$ 
denote by $\sigma_A[{\cal I}]$ the smallest $\sigma$-algebra
in the measurable space $(A,{\cal F}|_A)$ containing ${\cal I}$.
With this notation we have $\sigma[{\cal I}]=\sigma_{\Omega}[{\cal I}]$
for ${\cal I}\subset {\cal F}$.
The following elementary property will be crucial
in the proof of the main result to come.
For a class of measurable set ${\cal I}\subset{\cal F}$ 
and a pair of measurable sets $A,B\in{\cal F}$ satisfying $B\subset A$,
denote the restriction of ${\cal I}$ on $A$ by
$\dsp {\cal I}|_A:=\{I\cap A\mid I\in{\cal I}\}$
and the insertion of $B$ into ${\cal I}_A$ by
\begin{equation}
\eqna{insertion}
{\cal I}_{A,B}
=\{ B\cap I \mid I\in{\cal I}|_A\}\cup\{ B\cup I \mid I\in{\cal I}|_A\}
\end{equation}
In particular, $\dsp {\cal I}_{A,A}={\cal I}|_A$ if $\Omega\in {\cal I}$.
\lemb
\lema{insertion}
Let $(\Omega,{\cal F})$ be a measurable space and assume that
a class of measurable sets ${\cal I}\subset{\cal F}$ satisfies 
\eqnu{sigmaalgebrawithtotallyorderedgeneratingclass}. Then
for any pair of measurable sets $A$ and $B$ satisfying $B\subset A$,
${\cal I}_{A,B}$ of \eqnu{insertion} satisfies 
\eqnu{sigmaalgebrawithtotallyorderedgeneratingclass} with
the total space $(\Omega,{\cal F})$ replaced by $(A,{\cal F}|_A)$.
\DDD
\leme
\prfb
Since ${\cal I}$ satisfies
\eqnu{sigmaalgebrawithtotallyorderedgeneratingclass},
it suffices to prove
$\dsp \sigma_A[{\cal I}_{A,B}]\supset {\cal F}|_A$\,,
all the other properties in
\eqnu{sigmaalgebrawithtotallyorderedgeneratingclass} with the 
substitution $(\Omega,{\cal F})=(A,{\cal F}|_A)$ being 
direct consequences of the assumptions.

Put $\dsp{\cal G}=\{F\in{\cal F}\mid F\cap A \in\sigma_A[{\cal I}_{A,B}]\}$.
Then as in a standard elementary argument in measure theory.
$\dsp \sigma_A[{\cal I}_{A,B}]\supset {\cal F}|_A$ is equivalent to
${\cal G}\supset {\cal F}$.
Since by assumption ${\cal F}=\sigma[{\cal I}]$,
it suffices to prove ${\cal G}\supset {\cal I}$ and that
it is a $\sigma$-algebra in $\Omega$.
The latter is a straightforward consequence of 
the definition of ${\cal G}$ and that
$\sigma_A[{\cal I}_{A,B}]$ is a $\sigma$-algebra in $A$.

Finally, to prove ${\cal G}\supset {\cal I}$, let $I\in {\cal I}$.
Then $B=B\cap A\in{\cal I}_{A,B}$ and 
$B\cap I=B\cap (I\cap A)\in{\cal I}_{A,B}$ imply
$B\cap I^c=B\cap (A\cap (B\cap I)^c)\in \sigma_A[{\cal I}_{A,B}]$,
hence, with $B\cup (I\cap A)\in{\cal I}_{A,B}$, it follows that
\[
A\cap I=(B\cup (I\cap A))\cap (A\cap (B\cap I^c)^c)
\in \sigma_A[{\cal I}_{A,B}],
\]
which implies $I\in {\cal G}$.
\QED
\prfe
We are ready to state the main theorem.
Note that as we defined above \eqnu{submodular},
we assume that $v$ is non-decreasing, continuous, and $v(\emptyset)=0$,
when we say that $v$ is a submodular or supermodular function.
\thmb
\thma{totorderimpliesgenDelbaenKusuoka}
Let $(\Omega,{\cal F})$ be a measurable space, and
assume that there exists ${\cal I}\subset {\cal F}$ satisfying
\eqnu{sigmaalgebrawithtotallyorderedgeneratingclass}.
Then for any submodular function $v:\ {\cal F}\to\reals$
and for any $A,B\in{\cal F}$ satisfying $B\subset A$,
there exists unique $\mu\in{\cal C}_{-,v}(A)$ such that
\begin{equation}
\eqna{ShapleyCore}
v(I)=\mu(I),\ I\in{\cal I}_{A,B}\,,
\end{equation}
where ${\cal C}_{-,v}(A)$ is as in \eqnu{submodularlocalcoredef}.
In particular, $v(B)=\mu(B)$ holds.

Consequently, \eqnu{submodulartotalcorecoherentriskdef} holds
for all $A,B\in{\cal F}$ satisfying $B\subset A$.
\DDD
\thme
We will prove \thmu{totorderimpliesgenDelbaenKusuoka} in \secu{proof}.

The corresponding result for supermodular function directly follows
from \thmu{totorderimpliesgenDelbaenKusuoka} and
\prpu{convexconcaveduality}.
\corb
\cora{totorderimpliesgenShapley}
Let $(\Omega,{\cal F})$ be a measurable space, and
assume that there exists ${\cal I}\subset {\cal F}$ satisfying
\eqnu{sigmaalgebrawithtotallyorderedgeneratingclass}.
Then for any supermodular function $v:\ {\cal F}\to\reals$
and for any $A,B\in{\cal F}$ satisfying $B\subset A$,
there exists unique $\mu\in{\cal C}_{+,v}(A)$ such that
\eqnu{ShapleyCore} holds,
where ${\cal C}_{+,v}(A)$ is as in \eqnu{localcoredef}.
In particular,  $v(B)=\mu(B)$ holds.

Consequently, \eqnu{totalcorecoherentriskdef} holds
for all $A,B\in{\cal F}$ satisfying $B\subset A$.
\DDD
\core
\thmu{totorderimpliesgenDelbaenKusuoka} and
\coru{totorderimpliesgenShapley} imply
corresponding results on Choquet integrable functions.
For a non-decreasing, continuous, and finite (real-valued)
set function $v:\ {\cal F}\to\reals$ on a measurable space
$(\Omega,{\cal F})$ and a measurable function $f:\ \Omega\to\reals$,
we define 
\begin{equation}
\eqna{Choquetgen}
v(f)=\lim_{y\to-\infty} \biggl(y\,v(\Omega)
+\int_{y}^{\infty} v(\{\omega\in\Omega\mid f(\omega)>z\})\,dz\biggr)
\end{equation}
whenever the right hand side is a real value
and we say that $f$ is $v$-integrable.
If either the Lebesgue integration or the limit diverge
in the right hand side of \eqnu{Choquetgen} we do not define $v(f)$.
If $v(f)$ of \eqnu{Choquetgen} is defined it is equal to 
the Choquet integration, 
or the asymmetric integral in terms of \cite[Chap.~5]{Denneberg}.
If in addition $v$ is submodular,
it is known \cite[Prop.~10.3]{Denneberg} that
\begin{equation}
\eqna{genDelbaenKusuoka}
v(f)=\sup_{\mu\in{\cal C}_{-,v}(\Omega)} \int_{\Omega} f\,d\mu
\end{equation}
holds.
(In the reference, the statements are for finitely additive mesures and
algeras, but the corresponding results hold for ($\sigma$-additive) measures
when working on measurable space with $\sigma$-algebra as we do here.)
The functional $\rho$ defined by $\rho(f)=v(-f)$ is called
the coherent risk measure in mathematical finance
and \eqnu{genDelbaenKusuoka} is known to be a basic formula 
\cite{risk99,risk01,riskmSK}.

The definition of ${\cal C}_{-,v}(\Omega)$ and monotonicity of
Choquet integration is known to imply
$\dsp v(f)\le\sup_{\mu\in{\cal C}_{-,v}(\Omega)} \int_{\Omega} f\,d\mu$.
Hence, in a similar spirit as in \thmu{totorderimpliesgenDelbaenKusuoka},
to prove \eqnu{genDelbaenKusuoka} it is sufficient to prove 
the existence of a measure $\mu\in{\cal C}_{-,v}(\Omega)$ such that
$\dsp v(f)=\int_{\Omega} f\,d\mu$ holds.
A proof in \cite{Denneberg} for the existence of such $\mu$
starts with considering 
a maximal set of comonotonic $v$-integrable functions ${\cal X}$.
With the aid of Hahn--Banach's Theorem,
the domain of the functional $v$ is then extended to
to the linear space ${\cal V}$ generated by ${\cal X}$,
as a linear non-negative functional $\tilde{v}$.
Then a set function $\mu$ defined by $\mu(A)=\tilde{v}(\chrfcn{A})$,
where $\chrfcn{A}$ denotes the indicator function of the set $A$,
is proved to be a measure.
Since the indicator functions of the level sets
$\dsp I_{f,z}=\{\omega\in\Omega\mid f(\omega)>z\}$, $z\in\reals$,
are comonotonic with $f$, $\dsp\mu(I_{f,z})=v(I_{f.z})$ holds,
so that the definition of Choquete integration imply
$\dsp v(f)=\int_{\Omega} f\,d\mu$.
We also note that this $\mu$ gives equality $v(B)=\mu(B)$ in
\thmu{totorderimpliesgenDelbaenKusuoka},
by choosing $f=\chrfcn{A}$.

We note that, since the above outlined proof in the reference uses
Hahn--Banach Theorem,
uniqueness of measure $\mu$ which attains the
equality in \eqnu{submodulartotalcorecoherentriskdef} does not follow
in general.
The uniqueness claims in
\thmu{totorderimpliesgenDelbaenKusuoka} and
\coru{totorderimpliesgenShapley} are consequences
of our proof in the next section which uses
the extension Theorem of measures from measures on
finitely additive algebra.

As a simple and direct application of 
\thmu{totorderimpliesgenDelbaenKusuoka},
we can state a follwing result for the Choquet integration $v(f)$ in
\eqnu{genDelbaenKusuoka}.
Denote the class of level sets by ${\cal I}_f=\{I_{f,z}\mid z\in\reals\}$.
Then \thmu{totorderimpliesgenDelbaenKusuoka} implies the following.
\corb
\cora{genDelbaenKusuoka}
Let $(\Omega,{\cal F})$ be a measurable space,
$v:\ {\cal F}\to\reals$ be a submodular function,
and $f:\ \Omega\to\reals$ be a $v$-integrable function.
If the class of level sets ${\cal I}_f$ satisfies 
\eqnu{sigmaalgebrawithtotallyorderedgeneratingclass}
then there exists a unique $\mu\in{\cal C}_{-,v}(\Omega)$ such that
$v(I_{f,z})=\mu(I_{f,z})$, $z\in\reals$.
Moreover, $\dsp v(f)=\int_{\Omega}f\,d\mu$ holds for this $\mu$,
hence \eqnu{genDelbaenKusuoka} also holds.
\DDD
\core

\section{Proof of main theorem}
\seca{proof}

\prfofb{\protect\thmu{totorderimpliesgenDelbaenKusuoka}}
Assume $A,B\in{\cal F}$ satisfy $A\subset B$.
We fix $A$ and $B$ throughout the proof.

Let ${\cal F}|_A$ be the restriction
\eqnu{restriction} to $A$ of ${\cal F}$, and
$\dsp{\cal I}_{A,B}$ be the insertion \eqnu{insertion} of $B$ 
to ${\cal I}|_A$\,.
\lemu{insertion} then implies that
$\emptyset,B,A\in{\cal I}_{A,B}$ and that
$\dsp {\cal I}_{A,B}$ is totally ordered with respect to inclusion
with $\dsp \sigma[{\cal I}_{A,B}]={\cal F}|_A$\,.

Denote by $\dsp {\cal J}_{A,B}$
the finitely additive class generated by $\dsp {\cal I}_{A,B}$\,.
Namely, $\dsp {\cal J}_{A,B}$ is an algebra of sets satisfying
$\dsp{\cal I}_{A,B}\subset {\cal J}_{A,B}\subset{\cal F}|_A$\,,
is closed under complement and union, and is the smallest class
with these properties.
Since $\dsp {\cal I}_{A,B}$ is totally ordered with respect to inclusion,
we have an explicit representation
\begin{equation}
\eqna{finiteadditiveclass}
{\cal J}_{A,B}=
\{ \bigcup_{i=1}^n (C_i\cap D_i^c) \mid
C_1\supset D_1\supset C_2\supset\cdots\supset D_n\,,
\ C_i,\ D_i\in{\cal I}_{A,B}\,,\ i=1,2,\ldots,n,\ \ n=1,2,3,\ldots
\}.
\end{equation}
Using the notation in the right hand side of \eqnu{finiteadditiveclass},
define a set function $\mu_{A,B}:\ {\cal J}_{A,B}\to\reals$ by
\begin{equation}
\eqna{finiteadditivemeasure}
\mu_{A,B}(\bigcup_{i=1}^n (C_i\cap D_i^c))=
\sum_{i=1}^n (v(C_i)-v(D_i)).
\end{equation}
The definition implies that $\mu_{A,B}$ is finitely additive.
We assume non-decreasing property and continuity in the definitions of
submodular and supermodular functions. 
Therefore $\mu_{A,B}$ is non-negative valued and
$\sigma$-additive on the finitely additive class $\dsp{\cal J}_{A,B}$\,.
The extension theorem of measures
implies that $\mu_{A,B}$ is uniquely extended to
a measure on 
$\dsp \sigma[{\cal J}_{A,B}]=\sigma[{\cal I}_{A,B}]={\cal F}|_A$.
We denote this measure by the same symbol so that
$\mu_{A,B}\in{\cal M}(A)$.
We complete a proof of the theorem by proving that
$\mu_{A,B}\in{\cal C}_{-,v}(A)$ and $v(B)=\mu_{A,B}(B)$.

Note that $\emptyset,B,A\in{\cal I}_{A,B}$\,.
We can therefore put 
$n=1$, $C_1=A$, $D_1=\emptyset$ in \eqnu{finiteadditivemeasure}
to obtain $\dsp\mu_{A,B}(A)=v(A)$.
Similarly, with $n=1$, $C_1=B$, $D_1=\emptyset$ 
in \eqnu{finiteadditivemeasure} we have $\dsp\mu_{A,B}(B)=v(B)$.
It remains to prove $\mu_{A,B}(E)\le v(E)$ for $\dsp E\in {\cal F}|_{A}$.

Let $\eps$ be an aribtrary positive real.
Since $E\in{\cal F}|_A$ implies $A\cap E^c\in{\cal F}|_A$\,,
the extension theorem of measures implies that there exists a 
countable union of sets in ${\cal J}_{A,B}$\,, 
which we denote by $K\in{\cal F}|_A$
such that
\begin{equation}
\eqna{totorderimpliesgenDelbaenKusuokaprf1}
A\cap E^c \subset K\subset A\ \mbox{ and }\ 
\mu_{A,B}(A\cap E^c)+\eps\ge \mu_{A,B}(K).
\end{equation}
We can use the expression in \eqnu{finiteadditiveclass} for sets in 
${\cal J}_{A,B}$ to express $K$ as
\begin{equation}
\eqna{totorderimpliesgenDelbaenKusuokaprf2}
K= \bigcup_{i\in\nintegers} (C_i\cap D_i^c),
\ C_i,\ D_i\in{\cal I}_{A,B}\,;\ C_i\supset D_i\,,\ i=1,2,\ldots.
\end{equation}
For each $n\in\nintegers$ put $\dsp K_n=\bigcup_{i=1}^n  (C_i\cap D_i^c)$.
Then
\begin{equation}
\eqna{totorderimpliesgenDelbaenKusuokaprf7}
K_n\in{\cal J}_{A,B}\,,\ n\in\nintegers,\ \ 
K_1\subset K_2\subset \cdots,\ \ \bigcup_{n\in\nintegers} K_n=K.
\end{equation}
$\dsp \mu_{A,B}(A\cap E^c)+\eps\ge \mu_{A,B}(K)$ in
\eqnu{totorderimpliesgenDelbaenKusuokaprf1} and 
\eqnu{totorderimpliesgenDelbaenKusuokaprf7} then imply
\begin{equation}
\eqna{totorderimpliesgenDelbaenKusuokaprf8}
\mu_{A,B}(E)-\eps\le \mu_{A,B}(A\cap K_n^c),\ n\in\nintegers,
\end{equation}
while $\dsp A\cap E^c \subset K\subset A$ in
\eqnu{totorderimpliesgenDelbaenKusuokaprf1}
and \eqnu{totorderimpliesgenDelbaenKusuokaprf7}
and monotonicity and continuity of $v$ imply
that there exists $n_0\in\nintegers$ such that
\begin{equation}
\eqna{totorderimpliesgenDelbaenKusuokaprf9}
v(E)+\eps\ge v(A\cap K_n^c),\ n=n_0,n_0+1.\ldots.
\end{equation}
Since $K_n\in{\cal J}_{A,B}$ and $\dsp{\cal J}_{A,B}$ is a finite algebra,
$A\cap K_n^c\in{\cal J}_{A,B}$, hence
we can use the expression \eqnu{finiteadditiveclass} and write
\begin{equation}
\eqna{totorderimpliesgenDelbaenKusuokaprf5}
A\cap K_n^c=\bigcup_{i=1}^{n'} (C'_i\cap (D'_i)^c);\ \ 
C'_1\supset D'_1\supset C'_2\supset\cdots\supset D'_{n'}\,,
\ C'_i,\ D'_i\in{\cal I}_{A,B}\,,\ i=1,2,\ldots,n'.
\end{equation}
$\dsp C'_i,\ D'_i\in{\cal I}_{A,B}$ implies
\begin{equation}
\eqna{totorderimpliesgenDelbaenKusuokaprf3}
\mu_{A,B}(C'_i)=v(C'_i),\ 
\mu_{A,B}(D'_i)=v(D'_i),\ 
i=1,2,\ldots.n'.
\end{equation}
Therefore,
\begin{equation}
\eqna{totorderimpliesgenDelbaenKusuokaprf6}
\mu_{A,B}(A\cap K_n^c)=\sum_{i=1}^{n'}\mu_{A,B}(C'_i\cap (D'_i)^c)
=\sum_{i=1}^{n'}(\mu_{A,B}(C'_i)-\mu_{A,B}(D'_i))
=\sum_{i=1}^{n'}(v(C'_i)-v(D'_i))
\end{equation}

Put 
\begin{equation}
\eqna{totorderimpliesgenDelbaenKusuokaprf4}
A_i=\bigcup_{j=i}^{n} (C'_j\cap (D'_{j})^c),\ \ i=1,2,\ldots,n,
\end{equation}
and $A_{n+1}=\emptyset$.
Then \eqnu{totorderimpliesgenDelbaenKusuokaprf5} and
\eqnu{totorderimpliesgenDelbaenKusuokaprf4} imply $A_1=A\cap K_n^c$\,,
and
\[
A_i\cup D'_i=C'_{i}\,,\ A_i\cap D'_i=A_{i+1}\,,\ \ i=1,2,\ldots,n,
\]
Since $v$ is submodular, we can apply
\eqnu{submodular} with $A=A_i$ and $B=D'_i$ to find
\[
 v(C'_{i})+v(A_{i+1})\le v(A_i)+v(D'_i),\ i=1,2,\ldots,n.
\]
Summing this up with $i$ and using
\eqnu{totorderimpliesgenDelbaenKusuokaprf6} leads to
$\dsp\mu_{A,B}(A\cap K_n^c)\le v(A\cap K_n^c)$.
This with 
\eqnu{totorderimpliesgenDelbaenKusuokaprf8} and
\eqnu{totorderimpliesgenDelbaenKusuokaprf9} implies
$\dsp\mu_{A,B}(E)-\eps\le v(E)+\eps$.
Since $\eps$ is an arbitrary positive real this implies
$\dsp\mu_{A,B}(E)\le v(E)$, which completes the proof.
\QED
\prfofe

\section{Example}
\seca{salggeneratingchain}

In this section we give examples of measure spaces $(\Omega,{\cal F})$ and 
${\cal I}$ satisfying \eqnu{sigmaalgebrawithtotallyorderedgeneratingclass}
so that \thmu{totorderimpliesgenDelbaenKusuoka} is applicable.

\subsection{Countable set.}

Consider the case $\dsp (\Omega,{\cal F})=(\nintegers, 2^{\nintegers})$,
and put $\dsp {\cal I}
=\{\{1,2,\ldots,i\}\mid i\in\nintegers\}\cup\{\emptyset,\nintegers\}$.
Then ${\cal I}$ satisfies 
\eqnu{sigmaalgebrawithtotallyorderedgeneratingclass},
hence \thmu{totorderimpliesgenDelbaenKusuoka} is applicable.

For example, assume further that $m,n\in\nintegers$ satisfy $m<n$, and
$b_i\in\nintegers$, $i=1,2,\ldots,m$, and
$c_j\in\nintegers$, $j=1,2,\ldots,{n-m}$, satisfy
\[
b_1<b_2<\cdots<b_m\,,\ \ 
c_1<c_2<\cdots<c_{n-m}\,,\ \mbox{and}\ 
c_j\ne b_i\,,\ i=1,\ldots,m,\  j=1,\ldots,n-m.
\]
Put $B=\{b_1,\ldots,b_m\}$ and $A=B\cup\{c_1,\ldots,c_{n-m}\}$.
Then \eqnu{insertion} implies
\[
{\cal I}_{A,B}= \{\emptyset\}\cup
\{\{b_1,\ldots,b_i\}\mid i=1,2,\ldots,m\}\cup
\{B\cup\{c_1,\ldots,c_j\}\mid j=1,2,\ldots,n-m\}.
\]
Define a measure $\dsp\mu$ defined on ${\cal F}|_A$ by
\[\arrb{l}\dsp
\mu(\{b_1\})=v(\{b_1\}),\\
\mu(\{b_i\})=v(\{b_1,\ldots,b_i\})-v(\{b_1,\ldots,b_{i-1}\}),
\ i=2,\ldots,m,\\
\mu(\{c_1\})=v(A\cup\{c_1\})-v(A),\\
\mu(\{c_j\})=v(A\cup\{c_1,\ldots,c_j\})-v(A\cup\{c_1,\ldots,c_{j-1}\}),
\ j=2,\ldots,n-m.
\arre\]
Then \eqnu{ShapleyCore} holds and
\thmu{totorderimpliesgenDelbaenKusuoka} also implies $\mu\in{\cal C}_{-,v}(A)$.

This reproduces, with the correspondence
between submodular functions and supermodular functions (convex games)
of \prpu{convexconcaveduality},
a classical theory of cores of convex games by Shapley \cite{71shapl}.

\subsection{$1$-dimensional Borel $\sigma$-algebra}

Consider the case $\dsp (\Omega,{\cal F})=([0,1), {\cal B}([0,1)))$,
where $[0,1)$ is a unit interval and ${\cal B}([0,1))$ denotes
the $\sigma$-algebra generated by the open sets in $[0,1)$.

Put $\dsp {\cal I}=\{[0,x)\mid 0\le x\le 1\}$.
Then ${\cal I}$ satisfies 
\eqnu{sigmaalgebrawithtotallyorderedgeneratingclass},
hence \thmu{totorderimpliesgenDelbaenKusuoka} is applicable.

For $A,B\in{\cal B}([0,1)))$ satisfying $B\subset A$,
we have, from \eqnu{insertion},
\[
{\cal I}_{A,B}=
\{ B\cap [0,x)\mid 0\le x\le1\}\cup\{ B\cup (A\cap [0,x))\mid 0\le x\le 1\}.
\]
For each $I\in {\cal I}_{A,B}$ put $\mu(I)=v(I)$.
Then \lemu{insertion} and
the proof of \thmu{totorderimpliesgenDelbaenKusuoka} in \secu{proof}
imply that $\mu$ is uniquely extended to a measure on ${\cal B}([0,1))|_A$
and $\mu\in{\cal C}_{-,v}(A)$.

\subsection{Borel algebra on Polish space}

Let $\Omega$ be a separable, complete, metric space (a Polish space) and 
${\cal B}(\Omega)$ be the Borel $\sigma$-algebra,
the  $\sigma$-algebra generated by the open sets.
The collection ${\cal J}$ of all the open balls,
each with radius of positive rationals and
the center chosen from a fixed dense subset of $\Omega$,
generates the class of open sets,
i.e., ${\cal J}$ is a countable open basis,
hence $\sigma[{\cal J}]={\cal B}(\Omega)$ holds.
It turns out that we can construct a class ${\cal I}$ 
which satisfies 
\eqnu{sigmaalgebrawithtotallyorderedgeneratingclass}
on $(\Omega,{\cal F})=(\Omega,{\cal B}(\Omega))$
from ${\cal J}$,
so that \thmu{totorderimpliesgenDelbaenKusuoka} is applicable.

More generally, we have the following.
\prpb
\prpa{sigmaalgebraisomorphism}
Let $(\Omega,{\cal F})$ be a measurable space
and assume that there exists a class of coutable measurable sets
${\cal J}=\{J_1,J_2,\ldots\}\subset {\cal F}$ such that 
$\sigma[{\cal J}]={\cal F}$.

Define a ${\cal F}/{\cal B}([0,1))$ measurable function 
$f:\ \Omega\to[0,1)$ by
\begin{equation}
\eqna{sepcompletemeasurerealsembedding}
f=\sum_{n=1}^{\infty} 3^{-n}\chrfcn{J_n}\,,
\end{equation}
and define a class of sets ${\cal I}\subset {\cal F}$ by
$\dsp {\cal I}=\{ f^{-1}([0,a)) \mid 0\le a\le 1 \}$,
where $\dsp f^{-1}([0,a))=\{\omega\in\Omega\mid f(\omega)< a\}$.
Then ${\cal I}$ satisies 
\eqnu{sigmaalgebrawithtotallyorderedgeneratingclass}.
\DDD
\prpe
\prfb
Since $f$ is measurable, $\sigma[{\cal I}]\subset {\cal F}$.
Since ${\cal F}=\sigma[{\cal J}]$, to prove ${\cal F}\subset \sigma[{\cal I}]$
It suffices to prove $J_N\in\sigma[{\cal I}]$ for each $N\in\nintegers$.
Note that the right hand side of
\eqnu{sepcompletemeasurerealsembedding} has a form of ternary expansion,
because the indicator function $\chrfcn{J_n}$ takes values in $\{0,1\}$.
We can therefore write, for each $N\in\nintegers$,
\[\arrb{l}\dsp
J_N=f^{-1}(\bigcup_{(a_1,\ldots,a_{N-1})\in\{0,1\}^{N-1}}
[\sum_{n=1}^{N-1}a_n3^{-n}+3^{-N},\ \sum_{n=1}^{N-1}a_n3^{-n}+2\cdot 3^{-N}))
\\ \dsp \phantom{J_N}
=\bigcup_{(a_1,\ldots,a_{N-1})\in\{0,1\}^{N-1}}
f^{-1}(\ [0,\ \sum_{n=1}^{N-1}a_n3^{-n}+2\cdot 3^{-N})\ )\cap
f^{-1}(\ [0,\ \sum_{n=1}^{N-1}a_n3^{-n}+3^{-N})\ )^c
\\ \dsp \phantom{J_N}
\in \sigma[{\cal I}].
\arre\]
This proves $\sigma[{\cal I}]={\cal F}$.
The remaining properties stated in
\eqnu{sigmaalgebrawithtotallyorderedgeneratingclass} follows
from $\dsp f^{-1}([0,1))=\Omega$, $\dsp f^{-1}(\emptyset)=\emptyset$, and
$\dsp 0\le a<a'\le 1 \Rightarrow f^{-1}([0,a))\subset f^{-1}([0,a'))$.
\QED
\prfe
We note that $L^{\infty}([0,1))$ is not seperable, hence 
\prpu{sigmaalgebraisomorphism} is not applicable.

\end{document}